\documentclass[12t, russian]{article}
\usepackage[cp1251]{inputenc}
\usepackage[T2A]{fontenc}
\usepackage[english,russian]{babel}
\usepackage{amsmath, amssymb}
\usepackage{multirow}

\newtheorem{theorem}{Theorem}

\newtheorem{corollary}{Corollary}

\voffset=- 30mm
\begin{document}
\sloppy
\large

\

\begin{center}
\textbf{Hadamard  and Vandermonde  determinants \\
and \\
Bernoulli  -- Euler  -- Lagrange  -- Aitken  -- Nikiporets  type numerical method \\
for roots of polynomials}
\end{center}
\bigskip
\bigskip

\begin{center}
\textbf{ \large M. M. Chernyavskij (1), A. V. Lebedev (2), Yu. V. Trubnikov (1)}
\end{center}

\bigskip

\begin{center}
{\large ((1) Vitebsk State University,  Belarus, \	(2) Belarusian State University,  Belarus)}
\end{center}

\bigskip
\bigskip
\bigskip

\begin{abstract}

{  In the  article we develop Euler -- Lagrange method  and calculate all the roots of an arbitrary complex polynomial $P(z)$ on the base of calculation (similar to the Bernoulli -- Aitken -- Nikiporets methods) of the limits of ratios of Hadamard determinants  built by means of coefficients of expansions into Taylor and Laurent series of the function~$\frac{P'(z)}{P(z)}$.}
\end{abstract}
\bigskip
\bigskip

{\bf Keywords:} \emph{\small root of a polynomial, Hadamard determinant, Vandermonde determinant, Taylor series, \\
Laurent~series}

\bigskip
\bigskip


\bigskip
\bigskip
\bigskip

Methods of numerical solutions for roots of polynomials in the direction discussed in this article have a long and thoughtful history.

In   1728 D.~Bernoulli  \cite{Bern} described a method which
bears his name of numerical solution for the largest in modulus real root of a polynomial with real coefficients
$P(x) = a_0{x^n} + {a_1}{x^{n - 1}} + ... + {a_{n - 1}}x + {a_n}, \ \ a_0, a_n \neq 0$ (i.e.  $P(x)$ is a polynomial of degree  $n$ not vanishing at 0). In this method calculation of the root
reduces to the calculation of the limit  of the sequence $\frac{t_{m+1}}{t_m}$ of ratios of neighbouring in numbers solutions to the difference equation
\begin{equation}\label{e-0}
a_0 t_m + a_1 t_{m-1} + ... + a_n t_{m-n}=0 \ \ m= n, n+1, ... ,
\end{equation}
built by means of coefficients of the polynomial  $P(x)$ (for details see, for example,  \cite{Mcn-Pan} Ch. 10). D.~Bernoulli  did not give a justification of his method.  In 1748 году L.~Euler in his book  \cite{Eul} devoted Chapter 17 to the analysis of Bernoulli's type method for numerical calculation of the largest (minimal) in modulus real root of a polynomial $P(x)$ that does not possess multiple roots. L.~Euler used power series (he called them recurrent series) built
for the function  $\frac{1}{P(x)}$ and calculated the limits of ratios of neighbouring coefficients of these series. He observed (by examples) that in the situation when $P(x)$ possesses a pair of the largest in modulus complex conjugated roots the method may not work -- the limit in question may not exist. In  1798 J. L. Lagrange developing Euler's ideas in \cite{Lag} described the corresponding method of calculation of the largest (minimal) in modulus real root of a polynomial $P(x)$, possessing multiple roots. He used the series built for the function $\frac{P'(x)}{P(x)}$. In 1927 A.C. Aitken \cite{Ait} generalised Bernoulli's method for calculation of the products of ordered in modulus real roots of $P(x)$. He used the limits of ratios of determinants built from successive in numbers solutions to the difference equation \eqref{e-0} (for details see, for example,  \cite{Mcn-Pan} Ch. 10, where a review of other similar in spirit methods of calculation of the roots of polynomials with real coefficients is contained as well). In the articles by
V. I. Shmoylov and  D. I. Savchenko  \cite{Sh-Sav} and  by V. I. Shmoylov and  G. A. Kirichenko \cite{Sh-Kir} on the base of developed by V. I. Shmoylov \cite{Shm}  $r/\varphi$-algorithm of summation of (diverging) continued fractions Aitken's method is converted into calculation of Nikiporets' continued fractions   $N_i^{(n)} := N_i(a_0,...,a_n)$ (ratios of infinite ``determinants'', expressed in terms of coefficients of $P(x)$). Namely for their calculation
the $r/\varphi$-algorithm is exploited.

\smallskip

In the present article we develop Euler -- Lagrange method  and calculate all the roots of an arbitrary complex polynomial $P(z)$ on the base of calculation of the limits of ratios of Hadamard determinants (similar to the Bernoulli -- Aitken -- Nikiporets methods) built by means of coefficients of expansions into Taylor and Laurent series of the function $\frac{P'(z)}{P(z)}$.

The corresponding methods for calculation of the largest (minimal) in modulus root of $P(z)$  were obtained in \cite{T-Ch-18,T-Ch-21}.

\medskip

Let  $P(z) = a_0{z^n} + {a_1}{z^{n - 1}} + ... + {a_{n - 1}}z + {a_n}, \ \ a_0,a_1,...,a_n \in {\mathbb{C}}; \  a_0, a_n \neq 0$ be an arbitrary polynomial of degree~$n$ not vanishing at 0. Thus  \begin{equation}\label{e-1}
P(z) = a_0(z-z_1)^{m_1}\cdot ...\cdot (z-z_p)^{m_p},
\end{equation}
where  $m_1+{m_2}+\ldots +{m_p}=n$ is the sum of multiplicities of the roots $z_j$, and $z_i \neq z_j$ for $i \neq j$, and $z_j~\neq 0, \ j=1,\dots , p$.
Along with $P(z)$ we consider a rational function
\begin{equation}\label{e-2}
\frac{P'( z)}{P(z)}=\sum\limits_{j=1}^p{\frac{m_j}{z-z_j}} = \sum\limits_{k = 0}^\infty  {c_k}{z^k}.
\end{equation}
Here the right hand part is the expansion of $\frac{P'( z)}{P(z)}$ into the Taylor series in the neighbourhood of 0.

\smallskip

Note at once that by the contemporary means of computer mathematics (eg.,  Maple or Wolfram Mathematica) one can in an elementary way calculate any number of coefficients of this series for an arbitrary given polynomial $P(z)$.

\smallskip

By the coefficients $c_k$ of the series \eqref{e-2} one can built Hadamard determinants. Namely,
for each pair of natural numbers $(k,r), \ k\geq 0, r>0$ the Hadamard determinant $H_{k,r}$ is given by
\begin{equation}\label{e-3}
 H_{k,r}:=
\left| \begin{matrix}
	{c_k} & {c_{k+1}}&\dots & {c_{k+r-1}}  \\
	{c_{k+1}} & {c_{k+2}} & \dots &{c_{k+r}}  \\
\dots &\dots &\dots &\dots \\
	{c_{k+r-1}} & {c_{k+r}} & \dots &{c_{k+2(r-1)}}  \\
\end{matrix} \right| \, .
\end{equation}

For a collection of numbers $(\alpha_1, \dots , \alpha_s), \ s>1$ the Vandermonde determinant  $V(\alpha_1, \dots , \alpha_s)$ is given by
\begin{equation}\label{e-4}
 V(\alpha_1, \dots , \alpha_s) :=
\left| \begin{matrix}
	{1} & {1}&\dots & {1}  \\
	{\alpha_1} & {\alpha_{2}} & \dots &{\alpha_{s}}  \\
\alpha_1^2 & \alpha_{2}^2 & \dots &\alpha_{s}^2\\
\dots &\dots &\dots &\dots \\
	\alpha_1^{s-1} & \alpha_{2}^{s-1} & \dots &\alpha_{s}^{s-1}  \\
\end{matrix} \right| \, ;
\end{equation}
where we set $V(\alpha_1)=1$.

Recall that $V(\alpha_1, \dots , \alpha_s) \neq 0$ iff $\alpha_i \neq \alpha_j$ for $i\neq j$.

By $\overline{V}(\alpha_1, \dots , \alpha_s)$ we denote the ``inversed'' Vandermonde determinant
\begin{equation}\label{e-4-0}
 \overline{V}(\alpha_1, \dots , \alpha_s) :=
\left| \begin{matrix}
\alpha_1^{s-1} & \alpha_{2}^{s-1} & \dots &\alpha_{s}^{s-1}  \\	
\dots &\dots &\dots &\dots \\
\alpha_1^2 & \alpha_{2}^2 & \dots &\alpha_{s}^2\\
	{\alpha_1} & {\alpha_{2}} & \dots &{\alpha_{s}}  \\	
{1} & {1}&\dots & {1}
\end{matrix} \right| \, ;
\end{equation}
and also set $\overline{V}(\alpha_1)=1$.

The properties of determinants imply the following relations between $V(\alpha_1, \dots , \alpha_s)$  and $\overline{V}(\alpha_1, \dots , \alpha_s)$:
\begin{equation}\label{e-4-1}
   V(\alpha_1, \dots , \alpha_s)   = (-1)^{\left[\frac{s}{2} \right]}\overline{V}(\alpha_1, \dots , \alpha_s) ,
\end{equation}
where $[x]$ is the integral part of the number $x$.  And if $\alpha_i \neq 0, \ i=1, ..., s$; then
\begin{equation}\label{e-4-2}
   V(\alpha_1, \dots , \alpha_s) = (\alpha_1 \cdot ... \cdot  \alpha_s)^{s-1} \overline{V}(\alpha_1^{-1}, \dots , \alpha_s^{-1}) = (\alpha_1 \cdot ... \cdot  \alpha_s)^{s-1} (-1)^{\left[\frac{s}{2} \right]} {V}(\alpha_1^{-1}, \dots , \alpha_s^{-1}) .
\end{equation}

\smallskip

The next statement relates Hadamard and Vandermonde determinants for the polynomial $P(z)$ under consideration.

\begin{theorem} \label{t-1}
Let  $(z_1,\dots, z_p)$ be the roots of the polynomial $P(z)$ \eqref{e-1}  and  $ \sum\limits_{k = 0}^\infty  {c_k}{z^k}$ be the Taylor series \eqref{e-2}. For any pair $(k,r), \ k\geq 0, 0< r \leq p$ the following equality holds
  \begin{equation}\label{e-8}
H_{k,r}  =(-1)^r r!\sum_{\begin{smallmatrix}
		j_1 <j_2<\dots < j_r\\
		 1\leq j_r \leq p
\end{smallmatrix}}
\frac{m_{j_1}\cdot ... \cdot m_{j_r}}{(z_{j_1}\cdot ... \cdot z_{j_r})^{k+2r-1}}
[{V}(z_{{j}_1}, \dots , z_{{j}_r})]^2 .
  \end{equation}

In particular,
\begin{equation}\label{e-9}
 H_{k,p} = (-1)^p \,p!\, m_{1}\cdot ... \cdot m_{p} \left(\frac{1}{z_{1}\cdot ... \cdot z_{p}} \right)^{k+2p-1} [V(z_{1}, \dots , z_{p})]^2\,.
\end{equation}

For  $r>p \ \ H_{k,r} =0$.
 \end{theorem}
P\,r\,o\,o\,f. \ From  \eqref{e-2} by a routine calculation  one obtains
$c_k = -\sum_{j=1}^{p}\frac{m_j}{z_j^{k+1}}$, and therefore

\begin{equation}\label{e-0-1}
 H_{k,r} =(-1)^r\left|
\begin{matrix}
	\sum\limits_{j=1}^{p}\frac{m_j}{z_j^{k+1}} & \sum\limits_{j=1}^{p}\frac{m_j}{z_j^{k+2}} & ... & \sum\limits_{j=1}^{p}\frac{m_j}{z_j^{k+r}}\\
	\sum\limits_{j=1}^{p}\frac{m_j}{z_j^{k+2}} & ... & ... &\sum\limits_{j=1}^{p}\frac{m_j}{z_j^{k+r+1}}  \\
...& ...& ... & ... \\
\sum\limits_{j=1}^{p}\frac{m_j}{z_j^{k+r}} & ... & ...  & \sum\limits_{j=1}^{p}\frac{m_j}{z_j^{k+2r-1}}  \\
\end{matrix} \right|\,.
 \end{equation}
Exploiting the determinants properties and taking into account that a determinant possessing proportional columns (lines) is equal to zero one concludes that \eqref{e-0-1} implies
\begin{equation*}\label{e-0-2}
 H_{k,r} =(-1)^r\sum_{\begin{smallmatrix}
		j_1,j_2,\dots,j_r\\
		 1\leq j_s \leq p \\
j_i \neq j_s
\end{smallmatrix}}
 \left|
\begin{matrix}
	\frac{m_{j_1}}{z_{j_1}^{k+1}} & \frac{m_{j_2}}{z_{j_2}^{k+2}} & ... & \frac{m_{j_r}}{z_{j_r}^{k+r}}\\
	\frac{m_{j_1}}{z_{j_1}^{k+2}} & ... & ... & \frac{m_{j_r}}{z_{j_r}^{k+r+1}}  \\
...& ...& ... & ... \\
\frac{m_{j_1}}{z_{j_1}^{k+r}} & ... & ...  & \frac{m_{j_r}}{z_{j_r}^{k+2r-1}}  \\
\end{matrix} \right|\,
 \end{equation*}
 \begin{equation*}\label{e-0-2}
  =(-1)^r\sum_{\begin{smallmatrix}
		j_1,j_2,\dots,j_r\\
		 1\leq j_s \leq p \\
j_i \neq j_s
\end{smallmatrix}}
\frac{m_{j_1}\cdot ... \cdot m_{j_r}}{z_{j_1}^{k+r}\cdot ... \cdot z_{j_r}^{k+2r-1}}
 \left|
\begin{matrix}
	{z_{j_1}^{r-1}} & {z_{j_2}^{r-1}} & ... & {z_{j_r}^{r-1}}\\
	{z_{j_1}^{r-2}} & ... & ... & {z_{j_r}^{r-2}}  \\
...& ...& ... & ... \\
1 & ... & ...  & 1  \\
\end{matrix} \right|\,
 \end{equation*}
\begin{equation}\label{e-0-3}
  =(-1)^r\sum_{\begin{smallmatrix}
		j_1,j_2,\dots,j_r\\
		 1\leq j_s \leq p \\
j_i \neq j_s
\end{smallmatrix}}
\frac{m_{j_1}\cdot ... \cdot m_{j_r}}{z_{j_1}^{k+r}\cdot ... \cdot z_{j_r}^{k+2r-1}}
(-1)^{{\rm sign}(j_1, ... , j_r)} \overline{V}(z_{{\bar j}_1}, \dots , z_{{\bar j}_r}),
  \end{equation}
where  $({\bar j}_1, {\bar j}_2,\dots , {\bar j}_r)$ is the ordering of the collection of numbers  $(j_1,j_2,\dots,j_r)$: \ ${\bar j}_1< {\bar j}_2<\dots < {\bar j}_r$ and  ${\rm sign}(j_1, ... , j_r)$ is the corresponding evenness of the permutation $(j_1, ... , j_r)$.

Note that
\begin{equation}\label{e-0-4}
  \sum_{\begin{smallmatrix}
		{\text{all the permutations}}\\
 (1,2,\dots, r)
 \end{smallmatrix}}
\frac{1}{1\cdot z_{i_2}z_{i_3}^2\cdot ... \cdot z_{i_r}^{r-1}}
(-1)^{{\rm sign}(i_1, ... , i_r)} = {V}(z_1^{-1}, \dots , z_r^{-1}) .
  \end{equation}
Now from \eqref{e-0-3}, and taking into account \eqref{e-0-4}, and relations \eqref{e-4-1}
and  \eqref{e-4-2} one obtains

\begin{equation*}\label{e-0-5}
 H_{k,r}  =(-1)^r\sum_{\begin{smallmatrix}
		j_1,j_2,\dots,j_r\\
		 1\leq j_s \leq p \\
j_i \neq j_s
\end{smallmatrix}}
\frac{m_{j_1}\cdot ... \cdot m_{j_r}}{z_{j_1}^{k+r}\cdot ... \cdot z_{j_r}^{k+2r-1}}
(-1)^{{\rm sign}(j_1, ... , j_r)} \overline{V}(z_{{\bar j}_1}, \dots , z_{{\bar j}_r})
  \end{equation*}
  \begin{equation*}\label{e-0-6}
  =(-1)^r\sum_{\begin{smallmatrix}
		j_1,j_2,\dots,j_r\\
		 1\leq j_s \leq p \\
j_i \neq j_s
\end{smallmatrix}}
\frac{m_{j_1}\cdot ... \cdot m_{j_r}}{(z_{j_1}\cdot ... \cdot z_{j_r})^{k+r}}\cdot \frac{1}{1\cdot z_{j_2}z_{j_3}^2\cdot ... \cdot z_{j_r}^{r-1}}
(-1)^{{\rm sign}(j_1, ... , j_r)} \overline{V}(z_{{\bar j}_1}, \dots , z_{{\bar j}_r})
  \end{equation*}
  \begin{equation*}\label{e-0-7}
  =(-1)^r\sum_{\begin{smallmatrix}
		j_1 <j_2<\dots < j_r\\
		 1\leq j_r \leq p
\end{smallmatrix}}
r!\,\frac{m_{j_1}\cdot ... \cdot m_{j_r}}{(z_{j_1}\cdot ... \cdot z_{j_r})^{k+r}}
{V}(z_{j_1}^{-1}, \dots , z_{j_r}^{-1})
 \overline{V}(z_{{j}_1}, \dots , z_{{j}_r})
  \end{equation*}
 \begin{equation*}\label{e-0-8}
  =(-1)^r r!\sum_{\begin{smallmatrix}
		j_1 <j_2<\dots < j_r\\
		 1\leq j_r \leq p
\end{smallmatrix}}
\frac{m_{j_1}\cdot ... \cdot m_{j_r}}{(z_{j_1}\cdot ... \cdot z_{j_r})^{k+r}} \,
\frac{1}{(z_{j_1}\cdot ... \cdot z_{j_r})^{r-1}} \cdot (-1)^{\left[\frac{r}{2} \right]}
{V}(z_{{j}_1}, \dots , z_{{j}_r}) \cdot (-1)^{\left[\frac{r}{2} \right]}
{V}(z_{{j}_1}, \dots , z_{{j}_r})
  \end{equation*}
  \begin{equation*}\label{e-0-9}
  =(-1)^r r!\sum_{\begin{smallmatrix}
		j_1 <j_2<\dots < j_r\\
		 1\leq j_r \leq p
\end{smallmatrix}}
\frac{m_{j_1}\cdot ... \cdot m_{j_r}}{(z_{j_1}\cdot ... \cdot z_{j_r})^{k+2r-1}}
[{V}(z_{{j}_1}, \dots , z_{{j}_r})]^2 .
  \end{equation*}
The proof is complete.

\medskip

The formula  \eqref{e-9}  implies that
\begin{equation}\label{e-10}
\frac{H_{k,p}}{H_{k+1,p}}  = z_{1}\cdot ...\cdot z_{p} \,.
\end{equation}
And for $r<p$ we have the following observation.

 \begin{theorem} \label{t-2}
 Let  $0<|z_1| \leq |z_2| \leq ... \leq |z_r| < |z_{r+1}| \leq |z_{r+2}|\leq ...\leq |z_p|$ $($for $r=p-1$ the condition is written as  $0<|z_1| \leq |z_2| \leq ... \leq |z_{p-1}|< |z_p|$$)$. Then
\begin{equation}\label{e-7}
  \lim_{k\to
 \infty} \frac{H_{k,r}}{H_{k+1,r}} = z_1\cdot ... \cdot z_r \,.
 \end{equation}
And herewith
 \begin{equation}\label{e-7-}
  \left|\frac{H_{k,r}}{H_{k+1,r}} - z_1\cdot ... \cdot z_r \right| < C q^{k+2r-1},
\end{equation}
where
 $$
 0 < q = \frac{|z_r|}{|z_{r+1}|} < 1 ,
 $$
i.e. the sequence  \eqref{e-7}  converges as a geometric progression.

And once $k$ is such that  $q^{k+2r}D< \varepsilon < \frac{1}{2}$, where
\begin{equation}\label{e-t2-*}
D = \sum_{\begin{smallmatrix}
		j_1 <j_2<\dots < j_r\\
		 1\leq j_r \leq p \\
(j_1,j_2,..., j_r) \, \neq \,  (1,2, ...,r)
\end{smallmatrix}}
d_{{j_1} ... {j_r}}, \qquad \qquad d_{{j_1} ... {j_r}} = \frac{m_{j_1}\cdot ... \cdot m_{j_r}}{m_{1}\cdot ... \cdot m_{r}}\cdot \left[\frac{{V}(z_{{j}_1}, \dots , z_{{j}_r})}{{V}(z_{1}, \dots , z_{r})}\right]^2,
\end{equation}
one can take  $C=|z_1\cdot ... \cdot z_r|2D(1+2\varepsilon)$.
\end{theorem}

P\,r\,o\,o\,f. \ By means of  \eqref{e-8} one has

\begin{equation*}
\frac{H_{k,r}}{H_{k+1,r}} = \frac{\sum_{\begin{smallmatrix}
		j_1 <j_2<\dots < j_r\\
		 1\leq j_r \leq p
\end{smallmatrix}}
\frac{m_{j_1}\cdot ... \cdot m_{j_r}}{(z_{j_1}\cdot ... \cdot z_{j_r})^{k+2r-1}}
[{V}(z_{{j}_1}, \dots , z_{{j}_r})]^2 }
{\sum_{\begin{smallmatrix}
		j_1 <j_2<\dots < j_r\\
		 1\leq j_r \leq p
\end{smallmatrix}}
\frac{m_{j_1}\cdot ... \cdot m_{j_r}}{(z_{j_1}\cdot ... \cdot z_{j_r})^{k+2r}}
[{V}(z_{{j}_1}, \dots , z_{{j}_r})]^2}
  \end{equation*}

\begin{equation}\label{e-t2-1}
 = (z_1\cdot ... \cdot z_r)\frac{\left[ 1+\sum_{\begin{smallmatrix}
		j_1 <j_2<\dots < j_r\\
		 1\leq j_r \leq p \\
(j_1,j_2,..., j_r) \, \neq \, (1,2, ...,r)
\end{smallmatrix}}
d_{{j_1} ... {j_r}} \, q_{{j_1} ... {j_r}}^{k+2r-1}\right]}
{\left[1+\sum_{\begin{smallmatrix}
		j_1 <j_2<\dots < j_r\\
		 1\leq j_r \leq p \\
(j_1,j_2,..., j_r) \, \neq \, (1,2, ...,r)
\end{smallmatrix}}
d_{{j_1} ... {j_r}} \, q_{{j_1} ... {j_r}}^{k+2r}\right]}\, ,
\end{equation}
where
\begin{equation}\label{e-t2-0}
d_{{j_1} ... {j_r}} = \frac{m_{j_1}\cdot ... \cdot m_{j_r}}{m_{1}\cdot ... \cdot m_{r}}\cdot \left[\frac{{V}(z_{{j}_1}, \dots , z_{{j}_r})}{{V}(z_{1}, \dots , z_{r}) }\right]^2, \qquad \qquad
q_{{j_1} ... {j_r}} =\frac{z_1\cdot ... \cdot z_r}{z_{j_1}\cdot ... \cdot z_{j_r}}\, .
\end{equation}
The conditions of the theorem imply that for  $(j_1,j_2,..., j_r) \neq (1,2, ...,r)$ one has
\begin{equation}\label{e-t2-2}
0<|q_{{j_1} ... {j_r}}| \leq \frac{|z_r|}{|z_{r+1}|} =: q <1.
\end{equation}
This along with  \eqref{e-t2-1},  and  \eqref{e-t2-0} implies
\begin{equation*}
  \lim_{k\to
 \infty} \frac{H_{k,r}}{H_{k+1,r}} = z_1\cdot ... \cdot z_r \,,
 \end{equation*}
i.e.  \eqref{e-7} is true.

Now let us verify the estimate  \eqref{e-7-}.

Exploiting  \eqref{e-t2-1}, \eqref{e-t2-0}, and \eqref{e-t2-2} one has
 \begin{equation*}
  \left|\frac{H_{k,r}}{H_{k+1,r}} - z_1\cdot ... \cdot z_r \right| =
\end{equation*}
\begin{equation}\label{e-t2-3}
 = \left|(z_1\cdot ... \cdot z_r)\frac{\left[ 1+\sum_{\begin{smallmatrix}
		j_1 <j_2<\dots < j_r\\
		 1\leq j_r \leq p \\
(j_1,j_2,..., j_r) \, \neq \, (1,2, ...,r)
\end{smallmatrix}}
d_{{j_1} ... {j_r}} \, q_{{j_1} ... {j_r}}^{k+2r-1}\right]}
{\left[1+\sum_{\begin{smallmatrix}
		j_1 <j_2<\dots < j_r\\
		 1\leq j_r \leq p \\
(j_1,j_2,..., j_r) \, \neq \, (1,2, ...,r)
\end{smallmatrix}}
d_{{j_1} ... {j_r}} \, q_{{j_1} ... {j_r}}^{k+2r}\right]} - z_1\cdot ... \cdot z_r \right| \, .
\end{equation}
From  \eqref{e-t2-3}, relaxing for brevity of the record  the indexes under the summation sign  $\sum$\,, one obtains
\begin{equation*}
  \left|\frac{H_{k,r}}{H_{k+1,r}} - z_1\cdot ... \cdot z_r \right|
 = \left|z_1\cdot ... \cdot z_r\right|\left|\frac{\left[ \sum
d_{{j_1} ... {j_r}} \, q_{{j_1} ... {j_r}}^{k+2r-1} - \sum
d_{{j_1} ... {j_r}} \, q_{{j_1} ... {j_r}}^{k+2r}\right]}
{\left[1+\sum
d_{{j_1} ... {j_r}} \, q_{{j_1} ... {j_r}}^{k+2r}\right]} \right|
\end{equation*}

\begin{equation*}
 \leq \left|z_1\cdot ... \cdot z_r\right|\frac{ \sum
d_{{j_1} ... {j_r}} \, \left|q_{{j_1} ... {j_r}}\right|^{k+2r-1} \left|1- q_{{j_1} ... {j_r}}\right| }
{\left|1-\sum
d_{{j_1} ... {j_r}} \, \left|q_{{j_1} ... {j_r}}\right|^{k+2r}\right|}
\end{equation*}

\begin{equation*}
 \leq \left|z_1\cdot ... \cdot z_r\right|\frac{ 2\sum
d_{{j_1} ... {j_r}} \, q^{k+2r-1} }
{\left|1-\sum
d_{{j_1} ... {j_r}} \, q^{k+2r}\right|}
\end{equation*}

\begin{equation}\label{e-t2-4}
 = \left|z_1\cdot ... \cdot z_r\right|\left(\frac{ 2\sum
d_{{j_1} ... {j_r}} \,  }
{1- q^{k+2r} \sum
d_{{j_1} ... {j_r}} \,}\right) \, q^{k+2r-1}
 \leq C \, q^{k+2r-1}\, ,
\end{equation}
that proves  \eqref{e-7-}.

Clearly the denominator  $\left(1- q^{k+2r} \sum
d_{{j_1} ... {j_r}}\right)$ in the latter expression is positive for sufficiently large $k$.
Introducing the notation  $D:= \sum d_{{j_1} ... {j_r}}$  we conclude that once $q^{k+2r}D< \varepsilon< \frac{1}{2}$, then  $\frac{1}{1 -q^{k+2r}D} < 1+2\varepsilon$.  And therefore
$$
\left|z_1\cdot ... \cdot z_r\right|\left(\frac{ 2\sum
d_{{j_1} ... {j_r}} \,  }
{1- q^{k+2r} \sum
d_{{j_1} ... {j_r}} \,}\right)  < \left|z_1\cdot ... \cdot z_r\right| 2D (1+2\varepsilon),
$$
thus one can take the constant  $C$ in  \eqref{e-t2-4} to be $\left|z_1\cdot ... \cdot z_r\right| 2D (1+2\varepsilon)$. The proof of the theorem is complete.

\smallskip
Note that  $H_{k,1} = c_k$. Therefore for calculation of the minimal in modulus root
one obtains the following statement that constitutes (for polynomials with real coefficients and their real roots) the essence of L.~Euler's observation in Chapter 17 \cite{Eul}. Euler did not give an estimate of the speed of approximations.
\begin{corollary}\label{l-1}
Let  $(z_1,\dots, z_p)$ -- the roots of the polynomial $P(z)$ \eqref{e-1}, $0<|z_1| < |z_2| \leq ... \leq |z_p|$ and $ \sum\limits_{k = 0}^\infty  {c_k}{z^k}$ is the Taylor series  \eqref{e-2}.
Then
\begin{equation}\label{e-7''}
  \lim_{k\to
 \infty} \frac{c_{k}}{c_{k+1}} = z_1\,.
 \end{equation}
And herewith
 \begin{equation}\label{e-7''-}
  \left|\frac{c_{k}}{c_{k+1}} - z_1 \right| < C q^{k+1},
\end{equation}
where
 $$
 0 < q = \frac{|z_1|}{|z_{2}|} < 1 ,
 $$
i.e. the sequence  \eqref{e-7''} converges as a geometric progression.

And once $k$ is such that $q^{k+2}(n-1)< \frac{1}{2}$,
one can take $C=|z_1|4(n-1)$.
\end{corollary}
P\,r\,o\,o\,f.  One needs only to verify the final formula for the constant   $C$. It follows from the estimates for $C$  in the statement of Theorem~\ref{t-2}. Namely, in the situation under consideration the formula \eqref{e-t2-*} implies
$$
D= \sum_{j=2}^p \frac{m_j}{m_1}  \leq n-1,
$$
and by the statement of Theorem~\ref{t-2} one can take
$C=|z_1|2D(1+2\cdot\frac{1}{2}))= |z_1|4(n-1)$.

\medskip

In essence Theorem~\ref{t-2} describes not only sufficient but also necessary conditions of existence of the limits under consideration. Namely, the next observation holds.

\begin{theorem}\label{t-3}
Let $0<|z_1| \leq |z_2| \leq ... \leq |z_r| = |z_{r+1}| \leq |z_{r+2}|\leq ...\leq |z_p|$. Then there does not exist a limit  $\lim_{k\to
 \infty} \frac{H_{k,r}}{H_{k+1,r}}$.
\end{theorem}

P\,r\,o\,o\,f.   In view of \eqref{e-t2-1} the existence (nonexistence) of a limit of the sequence $\frac{H_{k,r}}{H_{k+1,r}}$ is equivalent to the existence (nonexistence) of a limit of the sequence

\begin{equation}\label{e-t3-1}
A_k := \frac{\left[ 1+\sum_{\begin{smallmatrix}
		j_1 <j_2<\dots < j_r\\
		 1\leq j_r \leq p \\
(j_1,j_2,..., j_r) \, \neq \, (1,2, ...,r)
\end{smallmatrix}}
d_{{j_1} ... {j_r}} \, q_{{j_1} ... {j_r}}^{k+2r-1}\right]}
{\left[1+\sum_{\begin{smallmatrix}
		j_1 <j_2<\dots < j_r\\
		 1\leq j_r \leq p \\
(j_1,j_2,..., j_r) \, \neq \, (1,2, ...,r)
\end{smallmatrix}}
d_{{j_1} ... {j_r}} \, q_{{j_1} ... {j_r}}^{k+2r}\right]}\, ,
\end{equation}
where  $d_{{j_1} ... {j_r}}$ and $ q_{{j_1} ... {j_r}}$ are described in  \eqref{e-t2-0}. The summands with  $|q_{{j_1} ... {j_r}}| <1$ do not influence the existence (nonexistence) of a limit of this sequence. By the condition of the theorem in the sums in
 \eqref{e-t3-1} there are summands with  $|q_{{j_1} ... {j_r}}| =1$, for example,
$|q_{{1} 2 ... {(r-1)}{(r+1)}}| =1$ and herewith $q_{{1} 2 ... {(r-1)}{(r+1)}}\neq 1$ since $z_r \neq z_{r+1}$.

Relaxing in  \eqref{e-t3-1} the summands with  $|q_{{j_1} ... {j_r}}| <1$, and denoting for brevity of the record multiindexes ${{j_1} ... {j_r}}$ by  $s$, one concludes that the existence of a limit of the sequence  $A_k$ \eqref{e-t3-1} is equivalent to the existence of a limit of the sequence
\begin{equation}\label{e-t3-2}
\tilde{A}_k := \frac{\left[ 1+\sum_s
d_s\, q_s^{k+2r-1}\right]}
{\left[1+\sum_s
d_s\, q_s^{k+2r}\right]}\, ,
\end{equation}
where  $|q_s|=1$ and there is $s_0$ such that $q_{s_0} \neq 1$.

Since  $|q_s|=1$ then  $q_s = e^{i\varphi_s}, \ 0<\varphi_s \leq 2\pi$.

One can come across the following two situations.

\smallskip

1) All  $\varphi_s$ are rationally commensurable with $2\pi$,  i.e.  $\frac{\varphi_s}{2\pi} = \frac{m_s}{n_s}, \ \ {m_s},{n_s}\in {\mathbb{N}}$.

In this case $\tilde{A}_k$ is a periodic sequence of period  $N= \text{LCM}\{n_s\}$ and it is not
a constant sequence as there is $s_0$ for which $q_{s_0} \neq 1$ (it can happen that some terms of this sequence are not defined, if  ${\left[1+\sum_s
d_s\, q_s^{k+2r}\right]} =0$). Thus there is no limit for $\tilde{A}_k$.

\smallskip

2) There is $\varphi_s$ which is rationally incommensurable with $2\pi$,  i.e.  $\frac{\varphi_s}{2\pi} \in {\mathbb R} \setminus {\mathbb Q}$.

Let us separate the indexes $s$ into two groups $\{s\}= \{t\}\sqcup \{v\}$, where  $\varphi_t$
are rationally commensurable with $2\pi$, and   $\varphi_v$ are rationally incommensurable with $2\pi$. With these notation
$\tilde{A}_k$ is written in the form
\begin{equation}\label{e-t3-4}
\tilde{A}_k := \frac{\left[ 1+\sum_t
d_t\, q_t^{k+2r-1} + \sum_v
d_v\, q_v^{k+2r-1}\right]}
{\left[1+\sum_t
d_t\, q_t^{k+2r} + \sum_v
d_v\, q_v^{k+2r}\right]}\, .
\end{equation}

Let $\frac{\varphi_t}{2\pi} = \frac{m_t}{n_t}, \ \ {m_t},{n_t}\in {\mathbb{N}}$ and  $N= \text{LCM}\{n_t\}$.

Consider the subsequence $\breve{A}_l := \tilde{A}_k, \ k +2r-1 = Nl, \ l= 1,2, ...$.
To finish the proof it is enough to establish nonexistence of a limit for $\breve{A}_l$.

By the choice of $N$ the sequence $\breve{A}_l $ has the form
\begin{equation}\label{e-t3-5}
 \breve{A}_l = \frac{\left[ C_1 + \sum_v
d_v\, q_v^{Nl}\right]}
{\left[C_2  + \sum_v
d_v\, q_v^{Nl+1}\right]}\, ,
\end{equation}
where  $C_1$, and $C_2$ are some constants.

Let $m$ be the number of indexes $v$,  and ${\bf T}^m$ be the $m$-dimensional torus in  $\mathbb{C}^m$: ${\bf T}^m = S^1\times ... \times S^1 = \{(\lambda_1, ..., \lambda_m): |\lambda_i| = 1, \ i=1,...,m \}$.
The collection $\{q_v^{N}\}_v$ is a point on the torus ${\bf T}^m$; and the closure of the set of the points $\{q_v^{Nl}\}_v, \ l=1,2, ...$ is a submanifold (isomorphic to a torus) of dimension  $m'\geq 1$ of the torus ${\bf T}^m$ ($m'$ is the number of rationally independent numbers in the collection $\{\frac{\varphi_v}{2\pi}\}_v$). This along with the explicit form \eqref{e-t3-5} of the sequence $\breve{A}_l$ implies nonexistence of a limit of this sequence. The proof is complete.

\smallskip

This theorem uncovers the noted in introduction L.~Euler's observation (\cite{Eul}, Ch.17) on the fact that under the existence (for a polynomial with real coefficients) of a pair of the largest in modulus complex conjugate roots the Bernoulli's type method may not work. Note herewith that the pairs of roots do not need to be complex conjugate (they can be anything -- and, in particular, real). As an example one can consider the polynomial $P(z) = z^2 -1$. Here
\begin{equation*}
\frac{P'( z)}{P(z)}={\frac{1}{z-1}} + {\frac{1}{z+1}}= \sum\limits_{k = 0}^\infty  {[(-1)^k -1]}{z^k}.
\end{equation*}
$H_{k,1}= [(-1)^k -1]$  and the sequence $\frac{H_{k,1}}{H_{k+1,1}}$ does not possess a limit.

\medskip

The results presented above give us a possibility to calculate the roots of a polynomial  $P(z)$
starting for the minimal in modulus  $0 < |z_1| < |z_2| < ...$.  Henceforth we describe the analogous procedure of calculation of the roots of a polynomial starting from the largest one.

Consider the expansion of the function $\frac{P'( z)}{P(z)}$  into the Laurent series in the neighbourhood of the infinity
(i.e. for $|z|> \max_{1\leq j \leq p}|z_j|$).
\begin{equation}\label{e-9'}
\frac{P'( z)}{P(z)}=\sum\limits_{j=1}^p{\frac{m_j}{z-z_j}} = \sum\limits_{k = 0}^\infty  \frac{b_k}{z^{k+1}}.
\end{equation}

For the coefficients of the series  \eqref{e-9'} one can built the corresponding Hadamar determinants. Namely, for each pair of natural numbers $(k,r), \ k\geq 0, r>0$ the Hadamar determinant  ${\bf{H}}_{k,r}$ is given by
\begin{equation}\label{e-3'}
 {\bf{H}}_{k,r}:=
\left| \begin{matrix}
	{b_k} & {b_{k+1}}&\dots & {b_{k+r-1}}  \\
	{b_{k+1}} & {b_{k+2}} & \dots &{b_{k+r}}  \\
\dots &\dots &\dots &\dots \\
	{b_{k+r-1}} & {b_{k+r}} & \dots &{b_{k+2(r-1)}}  \\
\end{matrix} \right| \, .
\end{equation}

An analogue of Theorem~\ref{t-1} for the Laurent series \eqref{e-9'} is the following

\begin{theorem} \label{t-4}
Let $(z_1,\dots, z_p)$ be the roots of a polynomial $P(z)$ \eqref{e-1}  and $ \sum\limits_{k = 0}^\infty  \frac{b_k}{z^{k+1}}$ be the Laurent series \eqref{e-9'}. For any pair $(k,r), \ k\geq 0, 0< r \leq p$ the following equality holds
 \begin{equation}\label{e-8'}
  {\bf{H}}_{k,r} = r!\sum_{\begin{smallmatrix}
		j_1<j_2<...<j_r\\
		 1\leq j_r \leq p
\end{smallmatrix}} m_{j_1}\cdot ... \cdot m_{j_r} \left({z_{j_1}\cdot ... \cdot z_{j_r}} \right)^{k} [V(z_{j_1}, \dots , z_{j_r})]^2\,.
 \end{equation}

In particular,
\begin{equation}\label{e-10}
 {\bf{H}}_{k,p} =  p!\, m_{1}\cdot ... \cdot m_{p} \left({z_{1}\cdot ... \cdot z_{p}} \right)^{k} [V(z_{1}, \dots , z_{p})]^2\,,
\end{equation}
For  $r>p \ \ {\bf{H}}_{k,r} =0$.
 \end{theorem}
P\,r\,o\,o\,f.  By an explicit computation one obtains from \eqref{e-9'}
that $b_k = \sum_{j=1}^{p}{m_j}{z_j^{k}}$, and therefore

\begin{equation}\label{e-t4-0-1}
 {\bf{H}}_{k,r} =\left|
\begin{matrix}
	\sum\limits_{j=1}^{p}{m_j}{z_j^{k}} & \sum\limits_{j=1}^{p}{m_j}{z_j^{k+1}} & ... & \sum\limits_{j=1}^{p}{m_j}{z_j^{k+r-1}}\\
	\sum\limits_{j=1}^{p}{m_j}{z_j^{k+1}} & ... & ... &\sum\limits_{j=1}^{p}{m_j}{z_j^{k+r}}  \\
...& ...& ... & ... \\
\sum\limits_{j=1}^{p}{m_j}{z_j^{k+r-1}} & ... & ...  & \sum\limits_{j=1}^{p}{m_j}{z_j^{k+2(r-1)}}  \\
\end{matrix} \right|\,.
 \end{equation}
Denoting $\xi_j := \frac{1}{z_j}$ one rewrites \eqref{e-t4-0-1} in the form
\begin{equation}\label{e-t4-0-2}
 {\bf{H}}_{k,r} =\left|
\begin{matrix}
	\sum\limits_{j=1}^{p}\frac{m_j}{\xi_j^{k}} & \sum\limits_{j=1}^{p}\frac{m_j}{\xi_j^{k+1}} & ... & \sum\limits_{j=1}^{p}\frac{m_j}{\xi_j^{k+r-1}}\\
	\sum\limits_{j=1}^{p}\frac{m_j}{\xi_j^{k+1}} & ... & ... &\sum\limits_{j=1}^{p}\frac{m_j}{\xi_j^{k+r}}  \\
...& ...& ... & ... \\
\sum\limits_{j=1}^{p}\frac{m_j}{\xi_j^{k+r-1}} & ... & ...  & \sum\limits_{j=1}^{p}\frac{m_j}{\xi_j^{k+2(r-1)}}  \\
\end{matrix} \right|\,.
 \end{equation}
Comparing  \eqref{e-t4-0-2}, and \eqref{e-0-1},  and using the formula \eqref{e-8} one concludes that

  \begin{equation*}\label{e-8''}
  {\bf{H}}_{k,r} = r!\sum_{\begin{smallmatrix}
		j_1<j_2<...<j_r\\
		 1\leq j_r \leq p
\end{smallmatrix}} m_{j_1}\cdot ... \cdot m_{j_r} \left({z_{j_1}\cdot ... \cdot z_{j_r}} \right)^{k+2(r-1)} [V(z_{j_1}^{-1}, \dots , z_{j_r}^{-1})]^2\,.
 \end{equation*}
This along with relations \eqref{e-4-2} between $V(z_{j_1}^{-1}, \dots , z_{j_r}^{-1})$ and
$V(z_{j_1}, \dots , z_{j_r})$ implies the equality \eqref{e-8'}. The proof is complete.

The formula  \eqref{e-10} implies that
\begin{equation}\label{e-10'}
\frac{{\bf{H}}_{k+1,p}}{{\bf{H}}_{k,p}}  = z_{1}\cdot ...\cdot z_{p} \,.
\end{equation}
And for $r<p$ one has the next analogue of Theorem~\ref{t-2}.

 \begin{theorem} \label{t-2'}
 Let  $|z_p| \geq |z_{p-1}| \geq ... \geq |z_{p-r+1}|> |z_{p-r}|  \geq |z_{p-r-1}|\geq ...\geq |z_1| >0$ $($for   $r=p-1$ the condition is written as  $0<|z_1| < |z_2| \leq ...  \leq |z_p|$$)$. Then
\begin{equation}\label{e-7'}
  \lim_{k\to
 \infty} \frac{{\bf{H}}_{k+1,r}}{{\bf{H}}_{k,r}} = z_{p-r+1}\cdot ... \cdot z_p \,.
 \end{equation}
And herewith
\begin{equation}\label{e-7'-}
  \left|\frac{{\bf{H}}_{k+1,r}}{{\bf{H}}_{k,r}} - z_{p-r+1}\cdot ... \cdot z_p \right| < C q^k,
\end{equation}
where
 $$
 0 < q  = \left|\frac{z_{p-r}}{z_{p-r+1}} \right| < 1 ,
 $$
 i.e. the sequence  \eqref{e-7'}  converges as a geometric progression.

And once $k$ is such that  $q^{k}D< \varepsilon < \frac{1}{2}$, where
\begin{equation}\label{e-t2-**}
D = \sum_{\begin{smallmatrix}
		j_1 >j_2> ...>  j_r\\
		 1\leq j_1 \leq p \\
j_1,j_2,..., j_r\, \neq \, p,p-1, ...,p-r+1
\end{smallmatrix}}
d_{{j_1} ... {j_r}}, \qquad \quad d_{{j_1} ... {j_r}} = \frac{m_{j_1}\cdot ... \cdot m_{j_r}}{m_{p}\cdot ... \cdot m_{p-r+1}}\cdot \left[\frac{{V}(z_{{j}_1}, \dots , z_{{j}_r})}{{V}(z_{p}, \dots , z_{p-r+1})}\right]^2,
\end{equation}
one can take  $C=|z_p\cdot ... \cdot z_{p-r+1}|2D(1+2\varepsilon)$.
\end{theorem}

P\,r\,o\,o\,f.  The proof goes along the scheme of the proof of Theorem~\ref{t-2}.

The formula  \eqref{e-8'} implies

\begin{equation*}
\frac{{\bf{H}}_{k+1,r}}{{\bf{H}}_{k,r}} = \frac{\sum_{\begin{smallmatrix}
		j_1 >j_2>... > j_r\\
		 1\leq j_1 \leq p
\end{smallmatrix}}
{m_{j_1}\cdot ... \cdot m_{j_r}}{(z_{j_1}\cdot ... \cdot z_{j_r})^{k+1}}
[{V}(z_{{j}_1}, \dots , z_{{j}_r})]^2 }
{\sum_{\begin{smallmatrix}
		j_1 >j_2>... > j_r\\
		 1\leq j_1 \leq p
\end{smallmatrix}}
{m_{j_1}\cdot ... \cdot m_{j_r}}{(z_{j_1}\cdot ... \cdot z_{j_r})^{k}}
[{V}(z_{{j}_1}, \dots , z_{{j}_r})]^2}
  \end{equation*}

\begin{equation}\label{e-t2-1'}
 = (z_p\cdot z_{p-1}\cdot ... \cdot z_{p-r+1})\frac{\left[ 1+\sum_{\begin{smallmatrix}
		j_1 >j_2> ... > j_r\\
		 1\leq j_1 \leq p \\
(j_1,j_2,..., j_r) \, \neq \, (p,p-1, ...,p-r+1)
\end{smallmatrix}}
d_{{j_1} ... {j_r}} \, q_{{j_1} ... {j_r}}^{k+1}\right]}
{\left[1+\sum_{\begin{smallmatrix}
		j_1 >j_2> ... > j_r\\
		 1\leq j_r \leq p \\
(j_1,j_2,..., j_r) \, \neq \, (p,p-1, ...,p-r+1)
\end{smallmatrix}}
d_{{j_1} ... {j_r}} \, q_{{j_1} ... {j_r}}^{k}\right]}\, ,
\end{equation}
where
\begin{equation}\label{e-t2-0'}
d_{{j_1} ... {j_r}} = \frac{m_{j_1}\cdot ... \cdot m_{j_r}}{m_{p}\cdot ... \cdot m_{p-r+1}}\cdot \left[\frac{{V}(z_{{j}_1}, \dots , z_{{j}_r})}{{V}(z_{p}, \dots , z_{p-r+1})}\right]^2, \qquad \quad
q_{{j_1} ... {j_r}} =\frac{z_{j_1}\cdot ... \cdot z_{j_r}}{z_p\cdot ... \cdot z_{p-r+1}}\, .
\end{equation}

From the conditions of the theorem it follows that for  $(j_1,j_2,..., j_r) \neq (p,p-1, ...,p-r+1)$ one has
\begin{equation}\label{e-t2-2'}
0<|q_{{j_1} ... {j_r}}| \leq \left|\frac{z_{p-r}}{z_{p-r+1}} \right| =: q <1.
\end{equation}
This along with  \eqref{e-t2-1'},  and \eqref{e-t2-0'} implies
\begin{equation*}
\lim_{k\to
 \infty} \frac{{\bf{H}}_{k+1,r}}{{\bf{H}}_{k,r}} = z_{p-r+1}\cdot ... \cdot z_p \,,
\end{equation*}
i.e. \eqref{e-7'} is true.

The estimate  \eqref{e-7'-} and the estimate for the constant $C$
is carried out according to the scheme of the proof of the estimate
\eqref{e-7-}. Namely, by the argument exploited in derivation of the estimate \eqref{e-t2-4},
and taking into account \eqref{e-t2-1'},   \eqref{e-t2-0'},    and  \eqref{e-t2-2'}, and relaxing for brevity of the record the indexes under the summation sign  $\sum$\,,  one obtains

\begin{equation*}
\left|\frac{{\bf{H}}_{k+1,r}}{{\bf{H}}_{k,r}} - z_{p-r+1}\cdot ... \cdot z_p \right|
 = \left|z_{p-r+1}\cdot ... \cdot z_p\right|\left|\frac{\left[ \sum
d_{{j_1} ... {j_r}} \, q_{{j_1} ... {j_r}}^{k+1} - \sum
d_{{j_1} ... {j_r}} \, q_{{j_1} ... {j_r}}^{k}\right]}
{\left[1+\sum
d_{{j_1} ... {j_r}} \, q_{{j_1} ... {j_r}}^{k}\right]} \right|
\end{equation*}

\begin{equation}\label{e-t2-4'}
 \leq \left|z_{p-r+1}\cdot ... \cdot z_p\right|\left(\frac{ 2\sum
d_{{j_1} ... {j_r}} \,  }
{1- q^{k} \sum
d_{{j_1} ... {j_r}} \,}\right) \, q^{k}
 \leq C \, q^{k}\, ,
\end{equation}
that proves  \eqref{e-7'-}.

Introducing the notation  $D:= \sum d_{{j_1} ... {j_r}}$  we conclude that once  $q^{k}D< \varepsilon< \frac{1}{2}$, then
$$
 \left|z_{p-r+1}\cdot ... \cdot z_p\right|\left(\frac{ 2\sum
d_{{j_1} ... {j_r}} \,  }
{1- q^{k} \sum
d_{{j_1} ... {j_r}} \,}\right) \, q^{k} < \left|z_{p-r+1}\cdot ... \cdot z_p\right|    2D (1+2\varepsilon)\,,
$$
that is one can take the constant  $C$ in  \eqref{e-t2-4'} to be equal $\left|z_{p-r+1}\cdot ... \cdot z_p\right|    2D (1+2\varepsilon)$. The proof is complete.

\medskip

Note that  ${\bf{H}}_{k,1} = b_k$. Therefore for the calculation of the largest in modulus
root one has the next (similar to Corollary~\ref{l-1}) statement that constitutes (for polynomials with real coefficients and their real roots) the essence of L.~Euler's observation in Chapter 17 \cite{Eul}. Euler did not give an estimate of the speed of approximations.

\begin{corollary}
Let $(z_1,\dots, z_p)$ be the roots of the polynomial $P(z)$ \eqref{e-1}, $|z_p| > |z_{p-1}| \geq ... \geq |z_1|>0$ and  $ \sum\limits_{k = 0}^\infty  \frac{b_k}{z^{k+1}}$ be the Laurent series  \eqref{e-9'}.
Then
\begin{equation}\label{e-7'''}
  \lim_{k\to
 \infty} \frac{b_{k+1}}{b_{k}} = z_p\,.
 \end{equation}
And herewith
 \begin{equation}\label{e-7''-}
  \left|\frac{b_{k+1}}{b_{k}} - z_p \right| < C q^{k},
\end{equation}
where
 $$
 0 < q = \frac{|z_{p-1}|}{|z_{p}|} < 1 ,
 $$
i.e. the sequence \eqref{e-7'''} converges as a geometric progression.

Once  $k$ is such that   $q^{k}(n-1)< \frac{1}{2}$,
one can take  $C=|z_p|4(n-1)$.
\end{corollary}
Here to derive the constant  $C$ we note that in the situation under consideration    \eqref{e-t2-**} implies
$$
D= \sum_{j=1}^{p-1} \frac{m_j}{m_p}  \leq n-1\,.
$$

\smallskip

Similar to Theorem~\ref{t-2}, Theorem~\ref{t-2'} in essence describes not only sufficient but also necessary conditions for existence of the limits under consideration. Namely, the next observation holds.

\begin{theorem}\label{t-3'}
Let $|z_p| \geq |z_{p-1}| \geq ... \geq |z_{p-r+1}| = |z_{p-r}|  \geq |z_{p-r-1}|\geq ...\geq |z_1| >0$. Then there does not exist a limit  $\lim_{k\to
 \infty} \frac{{\bf{H}}_{k+1,r}}{{\bf{H}}_{k,r}}$.
\end{theorem}

The proof can be derived by the same argument as the proof of Theorem~\ref{t-3}.

\end{document}